\documentclass[12pt]{amsart}
\usepackage{scrextend}
\usepackage[english]{babel}
\usepackage{amsmath}
\usepackage{amsthm}
\usepackage{amssymb}
\usepackage{mathrsfs}
\usepackage{mathdots}
\usepackage{enumerate,MnSymbol}
\usepackage[notcite, final, notref]{showkeys}
\usepackage{dsfont}
\usepackage{color}
\usepackage{shadow}
\newtheorem{theorem}{Theorem}[section]
\newtheorem{problem}[theorem]{Problem}
\newtheorem{lemma}[theorem]{Lemma}
\newtheorem{proposition}[theorem]{Proposition}
\newtheorem{corollary}[theorem]{Corollary}
\newtheorem{definition}[theorem]{Definition}

\theoremstyle{definition}
\newtheorem{remark}[theorem]{Remark}

\newcommand{\e}{\epsilon}

\newcommand{\ran}{{\rm Ran}\,}

\title[]{$q$-Rational functions and interpolation with complete Nevanlinna Pick Kernels}

\author[D. Alpay]{Daniel Alpay}
\address{(DA)
Faculty of Mathematics, Physics, and Computation\\
Schmid College of Science and Technology\\
Chapman University\\
One University Drive
Orange, California 92866\\
USA}
\email{alpay@chapman.edu}

\author[P. Cerejeiras]{Paula Cerejeiras}
\address{(PC) Center for research and development in mathematics and applications\\Department of mathematics, University of Aveiro  \\
Campus Universit\'ario de Santiago  \\ 3810-193 Aveiro\\Portugal}
\email{pceres@ua.pt}

\author[U. Kaehler]{Uwe Kaehler}
\address{(UK) Center for research and development in mathematics and applications\\Department of mathematics, University of Aveiro  \\
Campus Universit\'ario de Santiago  \\ 3810-193 Aveiro\\Portugal}
\email{ukaehler@ua.pt}

\author[B. Schneider]{Baruch Schneider}
\address{(BS) University of Ostrava\\ Department of Mathematics\\
30.dubna 22, 70200 Ostrava\\Czech Republic }
\email{baruch.schneider@osu.cz}

\begin{document}
\maketitle

\begin{abstract}
In this paper we introduce the concept of matrix-valued $q$-rational functions. In comparison to the classic case we give different characterizations with principal emphasise on realizations and discuss algebraic manipulations. We also study the concept of Schur multipliers and complete Nevanlinna Pick kernels in this context and provide first applications in terms of an interpolation problem using Schur multipliers and complete Nevanlinna Pick kernels. 
\end{abstract}
\tableofcontents

{\bf keywords:} rational functions; $q$-calculus; Jordan chains; CNP kernels; Schur multipliers.\\


\section{Introduction}
\setcounter{equation}{0}
The notion of rational function, Hardy space and backward-shift operator play an important role in classical analysis and its applications to domains
such as operator models, linear system theory, inverse scattering and others. This circle of ideas can be called Schur analysis, and
englobes (at least, and in a non-exhaustive list):\smallskip

{\sl The Hardy space,} and the underlying reproducing kernel. The Fock space plays also an increasingly important role in this field.

{\sl Shifts,} the forward-shift, backward-shift, Beurling -Lax theorem and generalizations of it.\\
{\sl  Schur multipliers,} and their co-isometric realizations.\\
  {\sl Rational functions,} various equivalent notions of rational function.\\
{\sl Interpolation problems,} for Schur multipliers and other related classes of functions.\\
{\sl System theory:} Links with linear system theory (filters, etc).\\
{\sl Inverse scattering:} and links with physics.\\

It is of interest to study the counterpart of these notions
in various domains such as several complex variables, hypercomplex analysis, calculus on diagonals, or linear systems indexed trees to cite a few.
In each case, one gets new notions and problems from the interactions with the field in question. We refer the reader to
\cite{MR2002b:47144,AJLV15,
arov2018multivariate,MR94b:47022,bgr,ball2021noncommutative,  Dym_CBMS, MR0255260,rosenbrock1970state} for background and further references.\smallskip

In our previous paper~\cite{Cerejeiras23} we began a study of Schur analysis in the setting of $q$-calculus. We defined in particular the counterpart of the Hardy space
and of the backward-shift operator  In the present paper we continue our study of Schur analysis in the $q$-calculus setting, and in particular define
and study Schur multipliers on the one hand and rational functions on the other.
Rational functions can be described in a number of equivalent ways, to be recalled in the sequel. We also study complete Nevanlinna Pick (CNP)
kernels (which are an important tool in Schur analysis) in this framework.\smallskip

Thus the present work is aimed to at least two audiences, focusing on the $q$-calculus  and on linear systems and related topics respectively. Researchers
working on the Fock space will also find results  of interest (see for instance the discussion pertaining to \eqref{rat3}).\\

In the setting of this introduction, it is convenient to recall the following characterization of a rational function (see Corollary \ref{cakb}):
A $\mathbb C^{m\times n}$-valued function analytic in a neighborhood of the origin, with Taylor series $f(z)=\sum_{k=0}^\infty z^kf_k$ is rational if and only if
there exists $N\in\mathbb N_0$ and  matrices $(C,A,B)\in \mathbb C^{m\times N} \times \mathbb C^{N\times N}\times \mathbb C^{N\times n}$ 
such that
\begin{equation}
  \label{rat1}
f_k=CA^{k}B,\quad k=0,1,2,\ldots
\end{equation}
(one could also have started with a formal power series; convergence in a neighborhood of the origin follows then from \eqref{rat1}).
We note that the power $A^k$ rather than $A^{k-1}$ is not a misprint in  \eqref{rat1} and that one has
\begin{equation}
  \sum_{k=0}^\infty z^kf_k=C(I_N-zA)^{-1}{B}.
\label{rat222}
\end{equation}
Let $q\in[0,1]$ and let $\left[k \right]_q$
\begin{equation}
  \left[k\right]_q= 1+q+\cdots +q^{k-1},\quad k>0,\quad{\rm and}\quad   [0]_q =1,
\label{kq111}
\end{equation}
and $ \left[k\right]_q!=\prod_{j=0}^k [j]_q$, i.e.
\begin{equation}
  \label{kq!}
 \left[k\right]_q!=1\cdot (1+q)\cdot (1+q+q^2)\cdots (1+q+\cdots +q^{k-1}).
\end{equation}
We will say that $f(z)$ is $q$-rational if now 
\begin{equation}
  \label{rat2}
  f_k=\frac{CA^{k}B}{ \left[k\right]_q!},\quad k=0,1,\ldots
\end{equation}
Note also that $A,B$ and $C$ are assumed independent of $q$. The case where they may depend on $q$ is also
discussed in the paper. When $q=0$ we get back to \eqref{rat1}; for $q=1$ the sequence \eqref{rat2} is the Borel transform (see Definition \ref{def5-4}) of the sequence
\eqref{rat1} and has been considered in  \cite[\S9]{daf1}.
Links with that work will be explained in the sequel.\smallskip

For $q\in(0,1)$ we will prove that
\begin{equation}
  \label{rat22}
  \sum_{k=0}^\infty z^kf_k=C\left(\prod_{j=0}^\infty (I_N-(1-q)zq^jA)^{-1}\right)B,
  \end{equation}
while for $q=1$
\begin{equation}
  \sum_{k=0}^\infty z^kf_k=Ce^{zA}B.
  \label{rat3}
\end{equation}
Expressions of the form \eqref{rat3} were considered in  \cite{daf1}.\\

When $n=m=N=1$ and $C=B=1$, $A=a\in\mathbb C$, formulas \eqref{rat222}, \eqref{rat22} and \eqref{rat3} reduce to
           \begin{equation}
             \label{eqza}
             f(z)=\begin{cases}\,\, \frac{1}{1-za},\quad \hspace{2.5cm}q=0,\\
               E_q(za)\,\, ({\rm see}\,\,\, \eqref{groningen}),\quad q\in(0,1),\\
              \, e^{za},\quad \hspace{2.75cm}q=1.
               \end{cases}
             \end{equation}

 The paper consists of 11 sections of which this introduction is the first.
 In Section \ref{sec2} we review some definitions and results on the $q$-calculus. It contains in particular the structure of the finite dimensional spaces invariant under the $q$-Jackson derivative; see \eqref{q-Jackson} for the latter and Proposition \ref{propo1234} for the structure result. After a short review of rational functions in Section 3 we are going to introduce $q$-rational functions in Sections 4 and 5 where we also discuss different characterizations. The question of Jordan chains as well as addition, multiplication, and inverses of rational functions will be  investigated in Sections 6 and 7. In the last sections we are introducing Schur multipliers and CNP kernels in the $q$-case with the question of interpolation by $q$-rational functions in mind.\\

Finally, a word on notation. ${\mathcal I}$ will denote the identity in a given underlying space (understood from the context). The identity
in $\mathbb C^{N\times N}$ will be denoted by $I_N$.
             
\section{Preliminaries on $q$-calculus}
\setcounter{equation}{0}
\label{sec2}
For general background on $q$-calculus we refer to \cite{Ernst2012,Kac2001}. In this section we review
some definitions and put them in the context of classical analysis (i.e. $q=0$). We also present a new result, Proposition \ref{propo1234},
on the counterpart of classic finite dimensional backwards-shift invariant spaces.\smallskip

The Hardy space $\mathbf H_2(\mathbb D)$ of the open unit disk $\mathbb D$ can be characterized (up to a multiplicative positive factor for the inner product) as the only Hilbert space of power series converging at the origin and such that $R_0^*=M_z$, where
$R_0$ is the backward-shift operator
\begin{equation}
  R_0f(z)=\begin{cases}\,\,\dfrac{f(z)-f(0)}{z},\quad z\not=0,\\
    \,\,f^\prime(0),\quad \hspace{1.2cm}z=0.\end{cases}
  \label{R0}
  \end{equation}
Note that in $\mathbf H_2(\mathbb D)$ we have the identities
  \begin{equation}
    R_0R_0^*={\mathcal I}
  \end{equation}
  and
  \begin{equation}
R_0M_z-M_zR_0={\mathcal I}-R_0^*R_0=C^*C
\end{equation}
where $Cf=f(0)$.\\

We here briefly review some of the results from our previous paper \cite{Cerejeiras23}.
For $q\in(0,1)$ we consider the $q$-Jackson derivative
  \begin{equation}\label{q-Jackson}
R_qf(z)=\frac{f(z)-f(qz)}{(1-q)z}.
\end{equation}
Furthermore
\begin{equation}
R_1=\partial.
  \end{equation}
  There, it is proved that the Hardy space $\mathbf H_{2,q}$ is to be the unique (up to a multiplicative positive constant) space of power series such
  that $R_q^*=M_z$. For $q=0$ and $q=1$ this formula reduces to the well-know formulas $R_0^*=M_z$ and $\partial^*=M_z$ respectively.\smallskip

We set:

 \begin{equation}
  E_q(z\overline{w})=\sum_{k=0}^\infty \frac{z^k\overline{w}^k}{ \left[k\right]_q!}.
  \label{groningen}
\end{equation}

\begin{definition}
  $\mathbf H_{2,q}$ is the reproducing kernel Hilbert space of functions analytic in $|z|<\frac{1}{1-q}$ with reproducing kernel
  $E_q(z\overline{w})$.
  \end{definition}

  When $q=0$ we get back the classical Hardy space of the open unit disk, with reproducing kernel $\frac{1}{1-z\overline{w}}$; see e.g. \cite{duren,MR1102893,MR924157}
  for the latter. When $q$ tends to $1$ we get back the Fock space, with reproducing kernel $e^{z\overline{w}}$.

\begin{lemma}
  $f(z)=\sum_{k=0}^\infty a_kz^k$ belongs to   $\mathbf H_{2,q}$ if and only if
  \begin{equation}
\sum_{k=0}^\infty \left[k\right]_q!\cdot |a_k|^2<\infty
    \end{equation}
  \end{lemma}

When $q=0$, the study of $R_0$-invariant subspaces of functions analytic at the origin play an important role in control theory, interpolation theory and related topics, especially in the case of matrix-valued functions. The case when the space is isometrically included in  $\mathbf H_2(\mathbb D)$ 
is of special interest.
  
\begin{proposition} Let $q\in[0,1)$.
  Let $\mathfrak M$ be a $N$-dimensional $R_q$-invariant subspace of $\mathbb C^n$-valued functions, analytic at the origin.
Let $F(z)$ be a $\mathbb C^{n\times N}$-valued function,with columns a basis of $\mathfrak M$. There exists a matrix
$A_q\in\mathbb C^{N\times N}$ (depending {\it a priori} on $q$) such that, with $C=F(0)\in\mathbb C^{n\times N}$
\begin{equation}
  \label{kobe}
  F(z)=C\prod_{j=0}^\infty (I_N-(1-q)zq^jA_q)^{-1}
\end{equation}
\label{propo1234}
\end{proposition}

\begin{proof}

\[
R_qF=FA_q,
\]
that is
\[
\frac{F(z)-F(qz)}{(1-q)z}=F(z)A_q
\]
or, equivalently
\[
F(z)(I_N-(1-q)zA_q)=F(qz)
\]
The result follows by iteration since the infinite product converges for $|q|<1$.
\end{proof}

The case of expressions \eqref{kobe} where $A_q$ does not depend on $q$ is of special interest.\\

 We conclude this section with:

\begin{lemma}\label{q-exponential}
  Let $A\in\mathbb C^{N\times N}$ with $\rho(A)<\frac{1}{1-q}$. Then,
  \begin{equation}
    \label{zwolle}
 \sum_{k=0}^\infty\frac{A^k}{ \left[k\right]_q!}=\prod_{j=0}^\infty(I_N-(1-q)q^jA)^{-1}
  \end{equation}
\end{lemma}

\begin{proof}
  For $A$ diagonalizable, this is a direct consequence of the scalar version \eqref{groningen} with $w=1$. But diagonalizable
  matrices are dense in $\mathbb C^{N\times N}$ and the result follows since both sides of \eqref{zwolle} are continuous
  functions of the entries of $A$ and that every square matrix is the limit of a sequence of diagonalizable matrices (see e.g. \cite[Exercise 3.1.28]{newbook}).
  \end{proof}

\section{Rational functions}
\setcounter{equation}{0}
In this section we recall various equivalent characterizations of rational functions in the classical setting.   It is well known
(see \cite{bgr,bgk1,MR569473,MR0255260}) that a $\mathbb C^{m \times n}$-valued function, say $F$, analytic in a neighborhood of
  the origin can be written as
  \begin{equation}
    \label{weiz}
    F(z)=D+zC(I_N-zA)^{-1}B,
\end{equation}
  where $N\in\mathbb N_0$ and  the matrices $(C,A,B)\in \mathbb C^{m\times N} \times \mathbb C^{N\times N}\times \mathbb C^{N\times n}$ 
  and $D=F(0)$.

     \begin{definition}
\label{rty1234}
       The pair $(C,A)\in\mathbb C^{m\times N}\times\mathbb C^{N\times N}$ is called observable if
       \begin{equation}
\bigcap_{\ell=0}^{N-1}\ker CA^\ell=\left\{0\right\}.
\end{equation}
       The pair $(A,B)\in\mathbb C^{N\times N}\times\mathbb C^{N\times n}$ is called controllable if
       \begin{equation}
         \bigcup_{\ell=0}^{N-1}\ran A^\ell B=\mathbb C^N.
       \end{equation}
      If the pair $(C,A)$ is observable and the pair $(A,B)$ is controllable, the triple $(C,A,B)\in\mathbb C^{m\times N}\times \mathbb C^{N\times N}\times \mathbb C^{N\times n}$ is called minimal.
         \end{definition}  
  
 A minimal realization is unique up to a  similarity matrix, i.e. the only degree of freedom is to replace $(C,A,B)$ by $(CT,T^{-1}AT, T^{-1}B)$ where $T$ is an
       arbitrary invertible function, called similarity matrix. Furthermore, any triple can be written in the form
         \begin{equation}
           \label{bgkkal}
A=\begin{pmatrix}* &* &*\\0&A_0&*\\ 0&0&*\end{pmatrix},\quad B=\begin{pmatrix}*\\ B_0\\ 0\end{pmatrix},\quad C=\begin{pmatrix}0&C_0&*\end{pmatrix},
           \end{equation}
           where the $*$ denote non-relevant entries and $(C_0,A_0,B_0)$ is minimal; furthermore,
           \[
             CA^nB=C_0A_0^nB_0,\quad n=0,1,\ldots,
             \]
             see e.g. \cite{bgk1}.\\

             We also mention that there  are other kind of realizations, centered at $\infty$, i.e.
             \[
             R(z)=D+C(zI_N-A)^{-1}B
\]
             or at a finite point, for instance. Here, together with realization \eqref{weiz} we will need another one, presented in the following lemma.

\begin{lemma}
  \label{vb}
  Let $R$ be a $\mathbb C^{m\times n}$-valued function analytic in a neighborhood of the origin. Then, $R$ is rational if and only if it can be written as
  \begin{equation}
    \label{weiz2}
    R(z)=C(I_N-zA)^{-1}B,
  \end{equation}
where $N\in\mathbb N_0$ and  the matrices $(C,A,B)\in \mathbb C^{m\times N} \times \mathbb C^{N\times N}\times \mathbb C^{N\times n}$.
  \end{lemma}

\begin{proof}
Apply \eqref{weiz} to 
  $zR(z)$ to write
  \[
    zR(z)=D+zC(I_N-zA)^{-1}B.
  \]
  Setting $z=0$ gives $D=0$. Dividing both sides of the equation by $z$ leads then to the result.
  \end{proof}

  \begin{corollary}
    Let $F$ be a $\mathbb C^{m\times n}$-valued function analytic in a neighborhood of the origin, with Taylor expansion $F(z)=\sum_{k=0}^\infty z^kF_k$.
    Then, $F$ is rational if and only if
    \begin{equation}
      \label{cakb}
      F_k=CA^kB,\quad k=0,1,\ldots
    \end{equation}
    where $(C, A,B)$ are as above.
  \end{corollary}

  \begin{proof}
    From \eqref{cakb} one obtains $\sum_{k=0}^\infty z^kR_k=C(I_N-zA)^{-1}B$, and the converge follows from Lemma \ref{vb} since, for $z$ small enough,
    $C(I_N-zA)^{-1}B=\sum_{k=0}^\infty z^kCA^kB$.
  \end{proof}

\begin{definition}
A $\mathbb C^{m\times N}$-valued function $T(z)$ analytic at the origin is called $R_0$-cyclic if the linear span $\mathcal M(T)$ of
the columns of the functions $R_0^\ell T$, $\ell=0,1,\ldots$ is finite dimensional.
\label{def56}
\end{definition}

\begin{proposition}\label{propo12345}
  $T$ is $R_0$-cyclic if and only if $T(z)=C(I_N-zA)^{-1}$
  where $(C,A)\in\mathbb C^{m\times N}\times\mathbb C^{N\times N}$ is observable.
  \end{proposition}

  For the proof, we refer to the literature, or to the proof of Theorem \ref{thr0} specialized to $q=0$.\\
  We now turn to the realizations of the product and the inverse of rational functions.

             \begin{proposition} (see \cite{bgk1})
               Let $F_j,$ $j=1,2$, be two matrix-valued functions analytic in a neighborhood of the origin, and respectively
               $\mathbb C^{n_1\times m}$-valued and $\mathbb C^{m\times m_2}$-valued. Let
               \begin{equation}
                 F_j(z)=D_j+zC_j(I_{N_j}-zA_j)^{-1}B_j,\quad j=1,2,
               \end{equation}
               be realizations of $F_1$ and $F_2$. Then $F = F_1 F_2$ has realization 
               \begin{equation}
                 \label{tg-2023}
F(z)=D+zC(I_N-zA)^{-1}B
               \end{equation}
               with
               \begin{equation}
                 \label{f1-fn}
                 \begin{split}
                   A&=\begin{pmatrix}A_1&B_1C_2\\ 0&A_2\end{pmatrix},\\
                   B&=\begin{pmatrix}B_1D_2\\B_2\end{pmatrix},\\
                   C&=\begin{pmatrix} C_1&  D_1C_2\end{pmatrix},\\
                   D&=D_1D_2.
                   \end{split}
                   \end{equation}
               \end{proposition}

 \begin{proposition} Let $F_j$, $j=1,2$, be two matrix-valued functions analytic in a neighborhood of the origin, and respectively
               $\mathbb C^{n_1\times m}$-valued and $\mathbb C^{m\times m_2}$-valued. Let
               \begin{equation}
                 F_j(z)=C_j(I_{N_j}-zA_j)^{-1}B_j,\quad j=1,2,
               \end{equation}
               be realizations of $F_1$ and $F_2$ of the form \eqref{weiz2}. Then $F = F_1 F_2$ has realization given by
           \begin{equation}
             F(z)=C(I_n-zA)^{-1}AB
           \end{equation}
           with
           \begin{eqnarray}
             A&=&\begin{pmatrix}A_1&B_1C_2\\ 0_{N_2\times N_1}&A_2\end{pmatrix}\\
                     B&=&\begin{pmatrix}0\\ B_2\end{pmatrix},\quad{\rm so ,\,\ that\,\,} AB=\begin{pmatrix}B_1(C_2B_2)\\ A_2B_2\end{pmatrix}\\
             C&=&\begin{pmatrix}C_1&0\end{pmatrix}.
                                     \end{eqnarray}
             \end{proposition}

             \begin{proof}
               We first remark that, with $A$ as above,
               \[
                 \begin{split}
                   \left(I_{N_1+N_2}-zA\right)^{-1}&=\begin{pmatrix}I_{N_1}-zA_1&-zB_1C_2\\ 0_{N_2\times N_1}&I_{N_2}-zA_2
                   \end{pmatrix}^{-1}\\
                  &= \begin{pmatrix}(I_{N_1}-zA_1)^{-1}&z(I_{N_1}-zA_1)^{-1}B_1C_2(I_{N_2}-zA_2)^{-1}\\
                     0_{N_1\times N_2}&(I_{N_2}-zA_2)^{-1}\end{pmatrix},
                   \end{split}
                 \]
                 and so
                 \[
C(I_N-zA)^{-1}B=z             C_1(I_{N_1}-zA_1)^{-1}B_1C_2(I_{N_2}-zA_2)^{-1}B_2.
\]
But $CB=0$ and so
\[
  \frac{1}{z}C(I_N-zA)^{-1}B=C(I_N-zA)^{-1}AB.
  \]
               \end{proof}


               These formulas can be iterated and will be used in the sequel.
               For a product of  $J$ matrix-functions $F_1,F_2,\ldots, F_J$ of compatible sizes, realization \eqref{weiz} becomes
                              \begin{equation}
                 \begin{split}
                   A&=\begin{pmatrix}A_1&B_1C_2 &B_1D_2C_3&\cdots &(B_1D_2\cdots D_{J-2}C_{J-2})&(B_1D_2\cdots D_{J-1}C_{J-1})
                     \\ 0&A_2&B_2C_3&\cdots &(B_2D_3\cdots D_{J-2}C_{J-2})&(B_2D_3\cdots D_{J-1}C_{J-1})\\
                     & &&\ddots           &          &                                                  &\\
                     & & & &A_{J-1}&B_{J-1}C_J\\
                     & & & &&A_J
                   \end{pmatrix},\\
                   B&=\begin{pmatrix}(B_1D_2\cdots D_J)\\(B_2D_3\cdots D_J)\\ \vdots\\ B_{J-1}D_J\\ B_J\end{pmatrix},\\
                   C&=\begin{pmatrix} C_1& D_1C_2&D_1D_2C_3&\cdots& (D_1\cdots D_{J-2}C_{J-1})&(D_1\cdots D_{J-1}C_J)\end{pmatrix},\\
                   D&=D_1D_2\cdots D_J
                   \end{split}
                   \end{equation}
               
                   while realization \eqref{weiz2} becomes:

                              \begin{equation}
                 \begin{split}
                   A&=\begin{pmatrix}A_1&(B_1C_2) &(B_1C_2)(B_2C_3)&\cdots &(B_1C_2)\cdots(B_{J-2}C_{J-1})&(B_1C_2)\cdots (B_{J-1}C_J)\\
                    0&A_2&A_2(B_2C_3)&\cdots &A_2(B_1C_2)\cdots(B_{J-2}C_{J-1})&A_2(B_2C_3)\cdots (B_{J-1}C_J)\\
                     & &&\ddots           &          &                                                  &\\
                     & & & &A_{J-1}&A_{J-1}(B_{J-1}C_J)\\
                     & & & &&A_J
                   \end{pmatrix},\\
                   B&=\begin{pmatrix}B_1(C_2B_2)(C_3B_3)\cdots (C_JB_J)\\
                     (A_2B_2)(C_3B_3)\cdots (C_JB_J)\\
                     \vdots\\
                     (A_JB_J)
                   \end{pmatrix},\\
                   C&=\begin{pmatrix} C_1& 0&0&0& \cdots&0\end{pmatrix}.
                   \end{split}
                   \end{equation}

                   Formulas \eqref{weiz} have been used in \cite{alpay2015realizations} in the theory of wavelet filters. We will use \eqref{weiz} and
                   \eqref{weiz2} in the study of $q$-rational functions. We first need to review the various characterizations of rational functions.  Before this we
                   present formulas for the inverse and recall 
\begin{equation}
F^{-1}(z)=D^{-1}-zD^{-1}C(I_N-z(A-BD^{-1}C))^{-1}BD^{-1}.
\label{inv-weiz}
\end{equation}
                        This formula is a direct consequence of the well-known formula for matrices of compatible sizes
                        \[
(I_u-MN)^{-1}I_u+M(I_v-NM)^{-1}M,\quad M\in \mathbb C^{u\times v},\,\,\,{\rm and},\,\,\,  N\in\mathbb C^{v\times u},
\]
applied to $M=zD^{-1}C$ and $N=(I_N-A)^{-1}B$. Here too we refer to \cite{bgk1} for further information and details. Finally we note that
formula \eqref{weiz2} above is not so convenient to compute the inverse. If $n=m$ and  $\det  CB\not =0$ we have using \eqref{inv-weiz}
  \[
    \begin{split}
      (C(I_N-zA)^{-1}B)^{-1}&=\left(CB+zC(I_N-zA)^{-1}AB\right)^{-1}\\
      &=(CB)^{-1}-z(CB)^{-1}\left(I_N-z(A-AB(CB)^{-1}C)\right)^{-1}AB(CB)^{-1}.
      \end{split}
    \]

In the classical case, a $\mathbb C^{n\times m}$-valued rational function $F(z)$ analytic in a neighborhood of the origin, with development
\[
  F(z)=\sum_{k=0}^\infty z^kF_k,
\]
can be defined in a number of equivalent ways:\\

{\bf 1. Polynomial characterization.}
A matrix-valued rational function can be written as
\[
F(z)=\frac{1}{p(z)}P(z)
\]
where $p$ is a scalar polynomial and $P$ is a $m\times n$ matrix with polynomial entries, and $p(0)\not=0$. \\

{\bf 2. Realization.} $F$ admits a realization, i.e. can be written as
\[
F(z)=C(I_N-zA)^{-1}B
\]
where $N\in\mathbb N_0$ and $(C, A,B)\in  \mathbb C^{m\times N} \times \mathbb C^{N\times N}\times \mathbb C^{N\times n}$. 
If $N=0$, $R(z)=CB$.\\

{\bf 3. Hankel operator.} The block Hankel operator
\[
  H=\begin{pmatrix}F_0 &F_1&F_2&\udots\\
    F_1&F_2&F_3&\udots\\
    F_2&F_3 &\udots &\udots\\
  \udots   &\udots & \udots  &\udots
  \end{pmatrix}
  \]
  has finite rank.\\

{\bf 4. Taylor series.}
The coefficients $F_k$ can be written as
\[
F_k=CA^{k}B,\quad k=0,1,2,\ldots
\]
for some matrices $(C, A,B)\in  \mathbb C^{m\times N} \times \mathbb C^{N\times N}\times \mathbb C^{N\times n}$. 
If $N=0$, all the $F_k=0$, $k=1,2,\ldots$.\\

{\bf 5. Backward-shift characterization.} With $R_0$ as in \eqref{R0}, the linear span of the columns of the functions $R_0^kF$, $k=0,1,2\ldots$ is finite dimensional.\\

The counterparts of these characterizations are given in Section \ref{sec5}.

 \section{$R_q$-cyclic functions}
\setcounter{equation}{0}
To define $q$-rational function we first need the counterparts of Definition \ref{def56} and Proposition \ref{propo12345}.

   \begin{definition} Let $q\in[0,1)$.
A $\mathbb C^{n\times m}$-valued function $T(z)$ analytic at the origin is called $R_q$-cyclic if the linear span $\mathcal M(T)$ of
       the columns of the functions $R_q^\ell T$, $\ell=0,1,\ldots$ is finite dimensional.
   \end{definition}

   \begin{theorem}
     \label{thr0}
     A  $\mathbb C^{n\times m}$-valued function $T(z)$ analytic at the origin is $R_q$-cyclic if and only if it can be written as
     \begin{equation}
       \label{montparnasse}
T(z)=C_q\prod_{j=0}^\infty (I_N-(1-q)zq^jA_q)^{-1}B_q
     \end{equation}
     for some matrices $(C_q, A_q,B_q)\in \mathbb C^{n\times N} \times \mathbb C^{N\times N}\times \mathbb C^{N\times m}$.
   \end{theorem}

   \begin{proof}
     Assume first that the linear span $\mathcal M(T)$ is finite dimensional, say of dimension $N$. Let $F$ be a $\mathbb C^{n\times N}$-valued function
     with columns a basis of $\mathcal M(T)$.      Thus there exists $B\in\mathbb C^{N\times m}$ such that $T(z)=F(z)B_q$.
     Since, by construction, $\mathcal M(T)$ is $R_q$-invariant, and hence $F$ is of the form \eqref{kobe}, and we get the result.\smallskip

     Conversely, assume $T$ to be of the form \eqref{montparnasse}. Then,
     \[
       \begin{split}
         (R_qT)(z)&=\frac{T(z)-T(qz)}{(1-q)z}\\
         &=C_q\left\{\frac{\prod_{j=0}^\infty (I_N-(1-q)zq^jA_q)^{-1}-\prod_{j=0}^\infty (I_N-(1-q)zq^{j+1}A_q)^{-1}}{(1-q)z}\right\}
         B_q\\
         &=C_q\left\{\frac{\prod_{j=0}^\infty (I_N-(1-q)zq^jA_q)^{-1}-\prod_{j=1}^\infty (I_N-(1-q)zq^jA_q)^{-1}}{(1-q)z}
         \right\}B_q\\
                  &=C_q\prod_{j=0}^\infty (I_N-(1-q)zq^jA_q)^{-1}\frac{I_N-(I_N-(1-q)zA_q)}{(1-q)z}B_q\\
                  &=C_q\prod_{j=0}^\infty (I_N-(1-q)zq^jA_q)^{-1}A_qB_q
       \end{split}
         \]
         and so
         \[
\mathcal M(T)\subset\left\{ F(z)= C_q\prod_{j=0}^\infty (I_N-(1-q)zq^jA_q)^{-1}c\,;\,c\in\cup_{\ell=0}^N\ran \, A_q^\ell B_q\right\}
\]
so that ${\rm dim}\, \mathcal M(T)<\infty$.
     \end{proof}

\section{$q$-rational functions}
\setcounter{equation}{0}
\label{sec5}
     \begin{definition} Let $q\in[0,1)$.
       A  $\mathbb C^{n\times m}$-valued function $T(z)$ analytic at the origin is called $q$-rational if it is
       $R_q$-cyclic with $(C_q,A_q,B_q)$ independent of $q$, i.e. of the form
       \begin{equation}
         \label{formula-rq-cyclic}
T(z,q)=C\prod_{j=0}^\infty (I_N-(1-q)zq^jA)^{-1}B.
\end{equation}
       \end{definition}

     A first example is given by $E_q(z)$ (see \eqref{groningen}).
     
       We note
\[
       \begin{split}
  T(z,0)&=C(I_N-zA)^{-1}B\\
  &=CB+zC(I_N-zA)^{-1}AB
\end{split}
\]
i.e. $T(z,0)$ is rational in the classical sense.\smallskip

The limiting case $q=1$ corresponds to expressions of the form \eqref{rat3}.\\

As stated before in the classic case we have different characterizations of rational functions.\\

Let us consider $F(z)=\sum_{k=0}^\infty \frac{z^k}{[k]_q!}F_k$. Then using Lemma~\ref{q-exponential} we have 
$$
F(z)=C\prod_{j=0}^\infty (I_N-(1-q)zq^jA)^{-1}B, 
$$
where $F_k=CA^kB$. Using formulas \eqref{weiz} and \eqref{weiz2}, and in the spirit of formulas given in \cite{alpay2015realizations}
in the setting of wavelets we have:

\begin{proposition}
  Let $F_J(z)=C\prod_{j=0}^J (I_N-(1-q)zq^jA)^{-1}B$.
  Then
  \[
    F_J(z)=D_J+zC_J(I_{(J+1)N}-zA_J)^{-1}B_J
  \]
  where
  \begin{eqnarray}
    A_J&=&(1-q)\begin{pmatrix} 1   & q&q^2&\cdots     &q^{J-1}&q^J\\                                                                                                                                 
            0& q&q^2&\cdots     &q^{J-1}&q^J\\
            0&0 & q^2 &\cdots     &q^{J-1}&q^J\\
                        \vdots&\vdots &\vdots & & \vdots & \vdots\\
            0& 0&0&\cdots     &q^{J-1}&q^J\\                                                                                                                                 
                        0& 0&0&\cdots     &0&q^J\end{pmatrix}\otimes A\\                                                                                                                                 
    B_J&=&\begin{pmatrix} B \\ B \\ \vdots \\ B\end{pmatrix}\quad (\mbox{$J+1$ copies of $I_N$})\\
C_J&=&(1-q)\begin{pmatrix}    1& q&q^2&\cdots     &q^{J-1}&q^J\end{pmatrix}\otimes CA\\
D_J&=&CB
    \end{eqnarray}
    \end{proposition}

  \begin{proof}
\[
(I_N-(1-q)zq^jA)^{-1}=I_N+z((1-q)q^jA)(I_N-z((1-q)q^jA))^{-1}
\]
Using formulas \eqref{f1-fn} (from $j=0$ to $j=J$ rather than from $j=1$ to $j=J$) we obtain
    \[
      \prod_{j=0}^J (I_N-(1-q)zq^jA)^{-1}=\mathscr D+z\mathscr C(I_{(J+1)N}-z\mathscr A)^{-1}\mathscr B,
      \]
      with
      \[
        \begin{split}
          \mathscr A&=(1-q)\begin{pmatrix} 1   & q&q^2&\cdots     &q^{J-1}&q^J\\                                                                                                                                 
            0& q&q^2&\cdots     &q^{J-1}&q^J\\
            0&0 & q^2&\cdots     &q^{J-1}&q^J\\
                        \vdots&\vdots &\vdots & & \vdots&\vdots\\
            0& 0&0&\cdots     &q^{J-1}&q^J\\                                                                                                                                 
                        0& 0&0&\cdots     &0&q^J\end{pmatrix}\otimes A\\                                                                                                                                 
          \mathscr B&=\begin{pmatrix} I_N\\ I_N\\ \vdots \\I_N\end{pmatrix}\quad (\mbox{$J+1$ copies of $I_N$})\\
          \mathscr C&=(1-q)\begin{pmatrix}    1& q&q^2&\cdots     &q^{J-1}&q^J\end{pmatrix}\otimes A\\
          \mathscr D&=I_p.
          \end{split}
        \]
        Note that $\mathscr C$ is the first block row of $\mathscr A$.
    \end{proof}

\begin{proposition} 
  Let $F_J(z)=C\prod_{j=0}^J (I_N-(1-q)zq^jA)^{-1}B$.
  Then
  \[
    F_J(z)=C_J(I_{(J+1)N}-zA_J)^{-1}B_J,
  \]
  where 
  \begin{eqnarray}
    A_J&=&\begin{pmatrix}(1-q)A&I_p&I_p&\cdots&I_p\\
      0&(1-q)qA&(1-q)qA   &\cdots&(1-q)qA\\
      & & & &\\
      \vdots& \vdots&\vdots & &(1-q)q^JA    \end{pmatrix},\\
B_J&=&\begin{pmatrix}B\\ (1-q)qAB\\ \vdots \\(1-q)q^JAB\end{pmatrix},\\
C_J&=&\begin{pmatrix} C&0&0&\cdots&0\end{pmatrix} \quad (\mbox{$J+1$ copies of $0_{p\times p}$}).
    \end{eqnarray}
\end{proposition}

    \begin{proof}
      Formulas \eqref{weiz2} specialized to \eqref{formula-rq-cyclic} give
\begin{equation}
  \begin{split}
    \mathscr A&=\begin{pmatrix}(1-q)A&I_p&I_p&\cdots&I_p\\
      0&(1-q)qA&(1-q)qA^\cdots &\cdots&(1-q)qA\\
      & & & &\\
      \vdots& & & &(1-q)q^JA    \end{pmatrix},\\
    \mathscr B&=\begin{pmatrix}I_p\\ (1-q)qA\\ \vdots \\(1-q)q^JA\end{pmatrix},\\
    \mathscr C&=\begin{pmatrix} I_p&0&0&\cdots&0\end{pmatrix}\quad (\mbox{$J+1$ copies of $0_{p\times p}$}).\\
  \end{split}
  \end{equation}
      \end{proof}

For studying the connection between classic rational functions and $q$-rational functions we need to establish the $q$-Borel transform. 

\begin{definition}
  \label{def5-4}
Let $F(z)=\sum_{k=0}^\infty z^kF_k$ be analytic in a neighbourhood of zero. We define the $q$-Borel transform of $F$ as $BF_q(z)=\sum_{k=0}^\infty \frac{z^k}{[k]_q!}F_k$.
\end{definition}

With the help of the $q$-Borel transform we can get a characterization of $q$-rational functions in terms of classic rational functions.

\begin{theorem}\label{q-Borel}
Let $F$ be an analytic function in the disk $|z|<R.$ Then its $q$-Borel transform is an analytic function in the disk $|z|<R/(1-q)$.
\end{theorem}

\begin{proof}
  Without loss of generality assume that for the coefficients of the series  $F(z)=\sum_{k=0}^\infty z^kF_k$ we have $\|F_k\|\neq 0$ for all $k$. This means that we have $\limsup_{k\rightarrow\infty}\frac{\|F_{k+1}\|}{\|F_k\|}=1/R$ so that for the series $BF_q(z)=\sum_{k=0}^\infty \frac{z^k}{[k]_q!}F_k$ we have
\[
  \limsup_{k\rightarrow\infty}\frac{\|F_{k+1}\|}{\|F_k\|[k+1]_q}=\limsup_{k\rightarrow\infty}\frac{\|F_{k+1}\|(1-q)}{\|F_k\|(1-q^{k+2})}=(1-q)/R.
  \]
\end{proof}

This theorem allows us to state the following corollary. 

\begin{corollary}
If $F$ is a function belonging to the Hardy space then its Borel transform $BF_q$ belongs to the $q$-Fock space $\mathbf{H}_{2,q}$.
\end{corollary}

\begin{proof}
  Let $F$ be a function in the Hardy space. Then its Borel transform $BF_q$ converges in the disk $|z|<1/(1-q)$. For its norm we have $\|BF_q\|^2_{\mathbf{H}_{2,q}}=
  \sum_{k=0}^\infty \|F_k/([k]_q!)\|^2 [k]_q!$. Since $F$ belongs to the Hardy space we have $\sum_{k=0}^\infty \|F_k\|^2<\infty$, which means that
  $\|BF_q\|_{\mathbf{H}_{2,q}}<\infty$.

\end{proof}

As a consequence we get the following characterization of a $q$-rational matrix-valued function.
\begin{theorem}\label{q-equiv}
A function $F$ is a $q$-rational matrix-valued function if and only if it is the $q$-Borel transform of a classic rational function. 
\end{theorem}

Since the Borel transform defines a one-to-one mapping the theorem is an immediate consequence of Theorem~\ref{q-Borel}.

Theorem~\ref{q-equiv} allows us to get similar characterizations of $q$-rational functions as in the case of rational functions. 

First of all, since we use the realization of a $q$-rational function as its definition by Lemma~\ref{q-exponential} we have that $T(z)=\sum_{k=0}^\infty z^k T_k$ with $T_k=\frac{CA^kB}{[k]_q!}$ (Taylor series characterization). From the fact that $R_qT=R_0BT$ we get that $\mathrm{span}\{R_q^kT, \,\, k=0,1,2, \ldots \}$
must be finite-dimensional as well as that the Hankel operator
\[
  H=\begin{pmatrix}{T_0}{[0]_q!} &{T_1}{[1]_q!}&{T_2}{[2]_q!}&\udots\\
   {T_1}{[1]_q!}& {T_2}{[2]_q!} & {T_3}{[3]_q!} &\udots\\
    {T_2}{[2]_q!}& {T_3}{[3]_q!} &\udots &\udots\\
  \udots  &\udots & \udots &\udots
  \end{pmatrix}
  \]
  has finite rank.
  Furthermore, since the corresponding characterizations are equivalent in the classic case they are also equivalent in the case of $q$-rational
  functions.\\

This leads to the following theorem.

\begin{theorem}
  For a power series $T(z)=\sum_{k=0}^\infty T_kz^k$  converging in a neighborhood of the origin, the following characterizations of $q$- rationality
  are equivalent:
\begin{enumerate}
\item {\bf Realization.} $T$ admits a realization, i.e. can be written as 
\[
T(z,q)=C\prod_{j=0}^\infty (I_N-(1-q)zq^jA)^{-1}B,
\]
with
where $N\in\mathbb N_0$ and $(C, A,B)\in \mathbb C^{n\times N} \times \mathbb C^{N\times N}\times \mathbb C^{N\times m}$. If $N=0$, $R(z)=CB$.\\

\item {\bf Hankel operator.} The block Hankel operator
\[
H=\begin{pmatrix}{T_0}{[0]_q!} &{T_1}{[1]_q!}&{T_2}{[2]_q!}&\udots\\
   {T_1}{[1]_q!}& {T_2}{[2]_q!} & {T_3}{[3]_q!} &\udots\\
    {T_2}{[2]_q!}& {T_3}{[3]_q!} &\udots &\udots\\
  \udots  &\udots & \udots &\udots
  \end{pmatrix}
  \]
  has finite rank.\\

\item {\bf Taylor series.}
The coefficients $T_k$ can be written as
\[
T_k=\frac{CA^{k}B}{[k]_q!},\quad k=0,1,2,\ldots
\]
for some matrices $(C, A,B)\in \mathbb C^{n\times N} \times \mathbb C^{N\times N}\times \mathbb C^{N\times m}$. 
If $N=0$, all the $T_k=0$, $k=1,2,\ldots$.\\

\item {\bf $q$-Jackson derivative characterization.} With $R_q$ as in \eqref{q-Jackson}, the linear span of the columns of the functions $R_q^kT$, $k=0,1,2\ldots$ is finite dimensional.\\

\end{enumerate}
\end{theorem}

\section{Jordan chains associated to $R_q$}
\setcounter{equation}{0}
In interpolation theory of classical Schur functions and  of related classes, Jordan chains in the interpolation data play a central role. To illustrate this, the trigonometric moment problem corresponds to $\left\{1,z,\ldots, z^N\right\}$, which is a Jordan chain associated to $R_0$. Missing one of the intermediate powers of $z$ will lead to a trigonometric moment with gaps, much more difficult to solve. In the present section we study the Jordan chains associated to $R_q,$ for $q\in(0,1)$. For $q=1$, Jordan chains correspond to the classical chains from the elementary theory of ordinary differential equations, for instance
\[
  \left\{e^x, xe^x,\ldots ,x^ne^x\right\}
  \]
is a Jordan chain associated to the eigenvalue $1$ of $\partial$ defined in $C^\infty(\mathbb R)$.

\begin{proposition}
Let $\mathfrak M$ be a $R_q$-invariant finite-dimensional space of dimension $N$, and let 
\[
F(z)=C_q\left(\prod_{j=0}^\infty (I_N-(1-q)zq^jA_q)^{-1}\right)=\begin{pmatrix}v_1(z)&v_2(z)&\cdots &v_N(z)\end{pmatrix}
\]
 be such that
\[
\mathfrak M=\left\{F(z)c\,\, ;\,\, c\in\mathbb C^N\right\}.
\]
Then $v_1,v_2,\ldots$ is a Jordan chain for $R_q$ if and only if $A_q$ is a Jordan cell.
\end{proposition}

\begin{proof}
Let $e_1,\ldots, e_N$ denote the canonical basis of $\mathbb C^N$, and assume first that $A_q$ is a $N\times N$ Jordan cell corresponding to the eigenvalue $\lambda$. We have:
\[
\begin{split}
A_qe_1&=\lambda e_1\\
A_qe_2&=\lambda e_2+e_1\\
&\vdots\\
A_qe_N&=\lambda e_N+e_{N-1}.
\end{split}
\]
Thus
\[
\begin{matrix}
R_qFe_1&=&FA_qe_1&=&\hspace{-1.9cm}\lambda Fe_1\\
R_qFe_2&=&FA_qe_2&=&\hspace{-7mm}\lambda Fe_2+Fe_1\\
&\vdots&& \vdots\\
R_qFe_N&=&FA_qe_N&=&\lambda Fe_N+Fe_{N-1}.
\end{matrix}
\]
Thus $v_1,\ldots, v_N$ is a Jordan cell for $R_q$.\smallskip

Conversely, assume that $\mathfrak M$ is spanned by the Jordan chain $v_1,\ldots, v_N$ based on $\lambda\in\mathbb C$, and let $f_1,\ldots, f_N\in\mathbb C^N$ be such that
\[
v_j=Ff_j,\quad j=1,\ldots, N.
\]
Thus,
\[
\begin{split}
R_qv_1&=\lambda v_1\\
R_qv_j&=\lambda v_1+v_{j-1},\quad j=2,\ldots N
\end{split}
\]
which ends the proof.
\end{proof}

As an example, consider the case $N=2$, and $\lambda=1$, so that
\[
A=  \begin{pmatrix}1&1\\0&1\end{pmatrix}=I_2+N,\quad N^2=0.
\]
Thus
\[
(I_2-wN)^{-1}=I_2+wN,\quad w\in\mathbb C,
  \]
and we have, with $C=\begin{pmatrix}1&1\end{pmatrix}$,
\[
  \begin{split}
    F(z)&=C\left(\prod_{j=0}^\infty\left(I_2-(1-q)zq^j\begin{pmatrix}1&1\\0&1\end{pmatrix}\right)^{-1}\right)\\
   &= C\left(\prod_{j=0}^\infty\left(I_2+(1-q)zq^j\begin{pmatrix}1&1\\0&1\end{pmatrix}\right)\right)\\
    &=
\left(\prod_{j=0}^\infty(1-(1-q)zq^j)\right)^{-1}C\prod_{j=0}^\infty \left(I_2+\frac{(1-q)zq^j}{1-(1-q)zq^j}N\right)\\
    &=\left(\prod_{j=0}^\infty(1-(1-q)zq^j)\right)^{-1}C\prod_{j=0}^\infty \begin{pmatrix}1&\sum_{j=0}^\infty\frac{(1-q)zq^j}{1-(1-q)zq^j}\\
      0&1\end{pmatrix}.
  \end{split}
\]
Hence,
\[
  \begin{split}
    F(z)\begin{pmatrix}1\\0\end{pmatrix}&=\frac{C\begin{pmatrix}1\\0\end{pmatrix}}{\prod_{j=0}^\infty(1-(1-q)zq^j)},\\
    F(z)\begin{pmatrix}0\\1\end{pmatrix}&=\frac{C\begin{pmatrix}
\sum_{j=0}^\infty\frac{(1-q)zq^j}{1-(1-q)zq^j}
        \\1\end{pmatrix}}{\prod_{j=0}^\infty(1-(1-q)zq^j)}\\
    &=\frac{\sum_{j=0}^\infty\frac{(1-q)zq^j}{1-(1-q)zq^j}}{{\prod_{j=0}^\infty(1-(1-q)zq^j)}}{C\begin{pmatrix}1\\0\end{pmatrix}}+
    \frac{1}{\prod_{j=0}^\infty(1-(1-q)zq^j)}C\begin{pmatrix}0\\1\end{pmatrix}.
  \end{split}
\]

When $q=0$, we get back to the classical formulas
\[
 \begin{split}
    F(z)\begin{pmatrix}1\\0\end{pmatrix}&=\frac{C\begin{pmatrix}1\\0\end{pmatrix}}{1-z},\\
    F(z)\begin{pmatrix}0\\1\end{pmatrix}&=\frac{zC\begin{pmatrix}   1    \\0\end{pmatrix}}{(1-z)^2}
    +\frac{C\begin{pmatrix}0\\1\end{pmatrix}}{1-z},
  \end{split}
 \]
for this example.

\begin{remark} One can also proceed as follows.
  Let

\begin{equation}
   A_q=\lambda I_N+ V,\quad V=\begin{pmatrix}0&1&0&\cdots &0\\
   0&0&1&\cdots &0\\
   \\
   \\
   0&0&\cdots &0&1\\
   0&0&\cdots&0&0\end{pmatrix},\quad V^N=0.
 \end{equation}

We write
\begin{equation}
 t_j=\frac{(1-q){\lambda}zq^j}{1-(1-q){\lambda}zq^j},
\end{equation}
and 
 \begin{eqnarray}
   b_1&=&\sum_{j=0}^\infty t_j\\
   b_2&=&\sum_{\substack{{j_1, j_2\ge 0}\\ j_1+j_2=2}}t_{j_1}t_{j_2}\\
   b_3&=&\sum_{\substack{{j_1,j_2,j_3\ge 0}\\ j_1+j_2+j_3=3}}t_{j_1}t_{j_2}t_{j_3}\\
   &\,\vdots \nonumber \\
   b_{N-1}&=& \sum_{\substack{{j_1,j_2,\cdots, j_{N-1}\ge 0}\\ j_1+j_2+\cdots+j_{N-1}=N-1}}t_{j_1}t_{j_2}\cdots t_{N-1}.
 \end{eqnarray}
then we have
\[
\begin{split}
  F(z)&=\left(\prod_{j=0}^\infty(1-(1-q){\lambda}zq^j)\right)^{-1}
  C\prod_{j=0}^\infty \left(I_N+\sum_{v=1}^{N-1}\left(\frac{(1-q){\lambda}zq^j}{1-(1-q){\lambda}zq^j}\right)^vV^v\right)\\
&=\left(\prod_{j=0}^\infty(1-(1-q){\lambda}zq^j)\right)^{-1}\left(\sum_{v=0}^{N-1}b_vV^v\right),
  \end{split}
 \]
\end{remark}


\section{Addition, multiplication and inverse}
\setcounter{equation}{0}
We now study addition, multiplication and inverse in the setting of $q$-rational functions. Since
\[
\frac{C_1A_1^kB_1}{ \left[k\right]_q!}+\frac{C_2A_2^kB_2}{ \left[k\right]_q!}=\frac{CA^kB}{ \left[k\right]_q!},
    \]
    with
    \[
      A=\begin{pmatrix}A_1&0\\ 0&A_2\end{pmatrix},\quad       B=\begin{pmatrix}B_1&B_2\end{pmatrix},\quad     C=\begin{pmatrix}C_1&C_2\end{pmatrix},\quad
    \]
    the sum of $q$ rational functions   is $q$-rational, but the product is not in general $q$-rational.

Let $R(z)$ be a matrix-valued rational function analytic at the origin,  with possibly non-minimal realization \eqref{weiz2}. We set
\begin{equation}
  T(R)(z,q)=C\left(\prod_{j=0}^\infty\left(I_N-z(1-q)q^jA\right)^{-1}\right)B.
  \label{rtyu}
\end{equation}
\begin{proposition}
   \eqref{rtyu} does not depend on the given realization.
               \end{proposition}

               \begin{proof}
                 In fact, using \eqref{bgkkal} one has
                 \[
                 T(R)(z,q)=C_0\left(\prod_{j=0}^\infty\left(I_N-z(1-q)q^jA_0\right)^{-1}\right)B_0.
               \]
               \end{proof}

               \begin{definition}
                 \label{defdef} We introduce the following multiplication (for rational functions of compatible sizes)
                 \begin{equation} 
                   T(R)(z,q)\star T(S)(z,q)=(T(R)(z,0)S(z,0))(z,q),
                 \end{equation}and its inverse (if it exists)
                 \begin{equation}
                   (T(R)^{-1})(z,q)=(T(R^{-1})(z,0)).
                 \end{equation}
                 \end{definition}

Addition is a special case of multiplication when one allows non-square factors since
\[
(R+S)(z)=\begin{pmatrix}R(z)&I_n\end{pmatrix}\begin{pmatrix}I_m\\ S(z)\end{pmatrix},
\]
where both $R$ and $S$ are $\mathbb C^{m\times n}$-valued.
Thus addition of $q$-rational functions can also be defined using Definition \ref{defdef}.

\section{Schur multipliers}
  \setcounter{equation}{0}
     \setcounter{equation}{0}
We can also study Schur multipliers in the context of $q$-rational functions. To this end 
we recall the well-known formula
\begin{equation}
  \label{blaschke-rkhs}
    \frac{1-b_a(z)\overline{b_a(w)}}{1-z\overline{w}}=\frac{1-|a|^2}{(1-z\overline{a})(1-\overline{w}a)}, \quad {\rm with}\quad
    b_a(z)=\frac{z-a}{1-z\overline{a}}.
\end{equation}
By the formula for the reproducing kernel in one dimension we can get the counterpart of \eqref{blaschke-rkhs} from
     \[
     \begin{split}
       K(z,w)&=\frac{E_q(za)E_q(\overline{wa})}{E_q|a|^2}\\
       &=\prod_{j=1}^\infty\frac{(1-(1-q)zq^ja)^{-1}(1-(1-q)\overline{w}q^j\overline{a})^{-1}}{(1-(1-q)q^j|a|^2)^{-1}}.
       \end{split}
            \]
     We write
     \[
       \frac{(1-(1-q)zq^ja)^{-1}(1-(1-q)\overline{w}q^j\overline{a})^{-1}}{1-(1-q)q^j|a|^2}=\frac{(1-uv)^{-1}(1-\overline{w}\overline{v})^{-1}}
       {(1-|v|^2)^{-1}},
       \]
       with
       \begin{eqnarray}
         u&=&\sqrt{1-q}q^{j/2}z,\\
                       w&=&\sqrt{1-q}q^{j/2}w,\\
         v&=&\sqrt{1-q}q^{j/2}a,
     \end{eqnarray}

and in particular we get for each $j$-term
\begin{equation}
  \begin{split}
        \frac{(1-(1-q)zq^ja)^{-1}(1-(1-q)\overline{w}q^j\overline{a})^{-1}}{(1-(1-q)q^j|a|^2)^{-1}}&=\frac{1-b_v(u)\overline{b_v(w)}}{1-u\overline{w}}\\
          &=\frac{1-b_{\sqrt{1-q}q^{j/2}a}(\sqrt{1-q}q^{j/2}z)\overline{b_{\sqrt{1-q}q^{j/2}a}(\sqrt{1-q}q^{j/2}w)}}{1-z\overline{w}(1-q)q^j}
          \label{4567890}
          \end{split}
\end{equation}
and so the one-dimensional reproducing kernel is
\begin{equation}
\frac{E_q(za)E_q(\overline{wa})}{E_q|a|^2}=\prod_{j=0}^\infty \frac{1-\left(b_{\sqrt{1-q}q^{j/2}a}(\sqrt{1-q}q^{j/2}z)\right)\left(\overline{b_{\sqrt{1-q}q^{j/2}a}(\sqrt{1-q}q^{j/2}w)}\right)}{1-z\overline{w}(1-q)q^j}
       \end{equation}

       Another way to extend \eqref{blaschke-rkhs} to the present setting is as follows.
  Recall that a (bounded) multiplier in $\mathbf H_{2,q}$ is a function, say $S$, analytic in $|z|<\frac{1}{1-q}$  and such that the operator $M_S$ of
  multiplication by $S$ is a bounded operator from $\mathbf H_{2,q}$ into itself. When the operator is a contraction, it will be called a
  Schur multiplier. When $q=1$, Liouville's theorem forces the set of bounded
  multipliers to be reduced to the zero function. On the other hand, when $q=0$, Schur multipliers are exactly functions analytic and contractive in
  the open unit disk, for which a vast theory originating with Schur's papers \cite{schur,schur2} exists.
  We now study the multipliers when $q\in[0,1)$. We first remark that formula \eqref{blaschke-rkhs} implies that

\begin{equation}
  \label{baab}
  \frac{1-b_a(\sqrt{1-q}    z)\overline{b_a(\sqrt{1-q}w)}}{1-(1-q)z\overline{w}}=\frac{1-|a|^2}{(1-\sqrt{1-q}z\overline{a})(1-
    \sqrt{1-q}\overline{w}a)}.
\end{equation}
 for $z,w\in\mathbb D_{1/\sqrt{1-q}}$. As a consequence we have:

\begin{proposition}
  Let $a\in\mathbb D$ and let $b_{a,q}(z)=b_a(\sqrt{1-q}z)$.
$M_{b_{a,q}}$ is a contractive multiplier of $\mathcal H(E_q)$.
\end{proposition}

  \begin{proof}
    Using \eqref{baab} we have:
    \[
      \begin{split}
        (1-b_{a,q}(z)\overline{b_{a,q}(w)})E_q(z\overline{w})&=
        \frac{(1-|a|^2)(1-(1-q)z\overline{w})}{(1-\sqrt{1-q}z\overline{a})
(1-\sqrt{1-q}\overline{w}a)}
          \frac{1}{\prod_{j=0}^\infty (1-z\overline{w}(1-q)q^j)}\\
        &= \underbrace{\frac{(1-|a|^2)}{(1-\sqrt{1-q}z\overline{a})
(1-\sqrt{1-q}\overline{w}a)}}_{{\rm positive\,\, definite}}
          \prod_{j=1}^\infty \underbrace{\frac{1}{1-z\overline{w}(1-q)q^j}}_{{\rm positive\,\, definite}}\\
      \end{split}
    \]
    which is a converging infinite product of positive functions which are positive definite in $|z|<\frac{1}{1-q}$, and therefore the kernel is positive definite there.
    \end{proof}

    The Blaschke factor $b_a$ is rational contractive in the open unit disk and unitay on the unit circle. A far reaching generalization of $b_a$
    which appears in the theory of interpolation of Schur functions is now reviewed, in order to present another family of Schur multipliers
    of $\mathbf H_{2,q}$. \smallskip

    We start with a signature matrix  $J\in\mathbb C^{N\times N}$  (i.e. satisfying $J=J^*=J^{-1}$) and an observable pair of matrices
    $(C,A)\in\mathbb C^{n\times N}\times\mathbb C^{N\times N}$ (see Definition \ref{rty1234}). Furthermore assume
    that the equation (called a Stein equation)
    \begin{equation}
P-A^*PA=C^*JC
      \end{equation}
      has a solution which is both Hermitian and invertible. Then, with $z_0$ a point in the resolvent set of $A$ and of modulus $1,$ we define
      \begin{equation}
\label{theta-gen}
        \Theta(z)=J-(1-z\overline{z_0})C(I_N-zA)^{-1}P^{-1}(I_N-z_0A)^{-*}C^*.
      \end{equation}

Consequently,       
      \begin{equation}
        \label{db-eq}
  \frac{J-\Theta(z)J\Theta(w)^*}{1-z\overline{w}}=C(I_N-zA)^{-1}P^{-1}(I_N-wA)^{-*}C^*,
\end{equation}
and so when $P>0$ the kernel \eqref{db-eq} is positive definite on the set $\Omega(\Theta)$
of complex $z$ which are either $z=0$ or such that $1/z\in\rho(A)$.
The function $\Theta$ is $J$-contractive in $\Omega(\Theta)\cap \mathbb D$ and $J$-unitary on the points of the unit circle where
it is defined.\smallskip

We note that the normalized Baschke factor $\frac{b_a(z)}{b_a(z_0)}$
and formula \eqref{blaschke-rkhs} correspond to \eqref{theta-gen} and \eqref{db-eq} with $A=\overline{a}$, $J=1$ and $C=1$.\smallskip

Using the same argument as for the Blaschke factor $b_a$ we obtain that the kernel
  \[
  \begin{split}
    \left(J-\Theta(\sqrt{1-q}z)J\Theta(\sqrt{1-q}w)^*\right)E_q(z\overline{w})&=\\
    &\hspace{-4cm}=C(I_N-\sqrt{1-q}zA)^{-1}P^{-1}(I_N-\sqrt{1-q}wA)^{-*}C^*          \frac{1}{\prod_{j=1}^\infty (1-z\overline{w}(1-q)q^j)}
  \end{split}
\]
is positive definite in $D_{1/\sqrt{1-q}}\cap\Omega(\Theta)$.\smallskip

By the characterization of the reproducing kernel Hilbert space associated to a product of functions positive definite on a common set
(see \cite[Theorem 2, page 361]{aron}) we can write:
\begin{theorem} We have

\[
  \mathcal H\left(\left(J-\Theta(\sqrt{1-q}z)J\Theta(\sqrt{1-q}w)^*\right)E_q(z\overline{w})\right)=\left\{f(z,z), \,\, f(z_1,z_2)\in H_1\otimes H_2
\right\}
  \]
where $H_1=\left\{f(\sqrt{1-q}z),\,\, f\in\mathcal H\left(\frac{J-\Theta(z)J\Theta(w)^*}{1-z\overline{w}}\right)\right\}$
  and $H_2$ is the reproducing kernel Hilbert
 space with reproducing kernel $\frac{1}{\prod_{j=1}^\infty (1-z\overline{w}(1-q)q^j)}$.
  \end{theorem}

After these examples we prove:
  
  \begin{theorem}
    Let $s$ be analytic in $|z|<\frac{1}{\sqrt{1-q}}$. Then $M_s$ is a contractive operator from $\mathbf H_{2,q}$ into itself if and only if
    $s(z\sqrt{1-q})$ is analytic and contractive in the open unit disk,i.e. if and only if $s$ is a Schur function.
  \end{theorem}

  \begin{proof}
    We have for every $z\in\mathbb D_{1/\sqrt{1-q}}=\left\{z\in\mathbb C,\,\,|z|<1/\sqrt{1-q}\right\}$,
    \[
(1-|s(z)|^2)E_q(|z|^2)\ge 0,\quad z\in\mathbb D_{1/\sqrt{1-q}}.
\]
i.e.,
\[
\left(  \frac{1-|s(z)|^2}{1-(1-q)|z|^2}\right)\cdot\frac{1}{\prod_{j=1}^\infty(1-(1-q)q^j|z|^2)}\ge 0,\quad z\in\mathbb D_{1/\sqrt{1-q}}.
\]
and in particular,
\[
 \frac{1-|s(z)|^2}{1-(1-q)|z|^2}\ge 0,\quad z\in\mathbb D_{1/\sqrt{1-q}}.
\]
Setting $\sigma(z)=s\left(\frac{z}{\sqrt{1-q}}\right)$, we have that the function $\sigma$ is analytic and contractive in the open unit disk, and
hence the result.
\end{proof}

\begin{corollary}
  Schur multipliers of $\mathbf H_{2,q}$ are the functions of the form $\sigma(\sqrt{1-q}z)$ where $\sigma$ is a Schur multiplier of
  $\mathbf H_2(\mathbb D)$.
\end{corollary}

     \begin{lemma}
Let $S$ be a Schur function. Then, for every $k$, the function $S(z\sqrt{1-q}q^{k/2})$ is a contractive multiplier of $\mathcal H(E_q)$.
\end{lemma}
\begin{proof}
This comes from the formula
\[
  \begin{split}
    (1-S(z\sqrt{1-q}q^{k/2})\overline{S(w\sqrt{1-q}q^{k/2})}E_q(z\overline{w})=&\\
    &\hspace{-3cm}=\frac{1-S(z\sqrt{1-q}q^{k/2})\overline{S(w\sqrt{1-q}q^{k/2})}}{1-z\overline{w}\sqrt{1-q}q^k}\frac{1}{\prod_{\substack{j=1\\j\not=k}}^\infty(1-z\overline{w}\sqrt{1-q}q^j)}
  \end{split}
\]
\end{proof}

Since a product of Schur multipliers is a Schur multiplier we have:
\begin{corollary}
  Let $S$ be a Schur function. Then, for every $0<n < m < \infty$ the product
  \[
\prod_{k=n}^{m}S(z\sqrt{1-q}q^{k/2})
  \]
  is a Schur multiplier of $\mathcal H(E_q)$.
  \end{corollary}


\section{Multipliers and CNP kernels}
\setcounter{equation}{0}
Let $K(z,w)$ be a positive definite kernel on the set $\Omega$, with associated reproducing kernel Hilbert space $\mathcal H(K)$, and let $S$ be a contractive multiplier (or
Schur multiplier), i.e. a (say $\mathbb C^{p\times q}$-valued) function $S$ defined on $\Omega$ and such that the operator of multiplication by $S$ is a contraction
from $\mathcal H(K)$ into itself. Equivalently, the function
\begin{equation}
  \label{ks}
  K_S(z,w)=(1-S(z)S(w)^*)K(z,w)
\end{equation}
is positive definite in $\Omega$. Let furthermore $(w_j,s_j)$, $j=1,\ldots ,m$ be pre-assigned in $\Omega\times\mathbb C^{p\times q}$.
Since \eqref{ks} is positive
definite, a necessary condition for a Schur multiplier to exist such that
\begin{equation}
  \label{inter}
  S(w_j)=S_j,\quad j=1,\ldots, m,
\end{equation}
is that the block-matrix
\begin{equation}
  \label{G}
G=\left((I_p-S_jS_\ell^*)K(w_j,w_\ell)\right)_{j,\ell=1}^m,
\end{equation}
is positive definite. This follows from
\[
M_S^* (K(\cdot, w)c)=K(\cdot, w)S(w)^*c.
  \]
This condition will not be necessary in general.

\begin{definition}
  Let $\Omega\subset\mathbb{C}$ be a domain. A function $K:\Omega\times\Omega\to\mathbb{C}$ is called a complete Nevanlinna-Pick kernel (or CNP-kernel) when
$G\ge 0$ is not only a necessary but also a sufficient condition for a Schur multiplier satisfying \eqref{inter} to exist.
\end{definition}

For the following result, see \cite{quiggin} and \cite{MR1882259}.

\begin{theorem}
  A positive definite kernel is a CNP-kernel if and only if $1/K(z,w)$ has one positive square. 
\end{theorem}

Thus a complete Nevanlinna-Pick kernel is of the form
\[
\frac{1}{a(z)\overline{a(w)}-\langle b(z),b(w)\rangle_{\ell_2(A,\mathbb C)}},
\]
where $A$ is an indexed set (in this work, $A=\left\{1,\ldots ,N\right\}$ or $\mathbb N$).
In the sequel we will only consider the case where $a(z)=1$.
This is no real loss of generality; if $a(z_0)=0$ for some $z_0\in\Omega$ then $b(z_0)=0$. So we consider kernels of the form
\begin{equation}
  \label{cnp00}
\frac{1}{1-\langle b(z),b(w)\rangle_{\ell_2(A,\mathbb C)}}.
\end{equation}

An important example for the sequel is
\[
\frac{1}{1-\langle z,w\rangle}
  \]
  where $z$ and $w$ run through the open unit ball of $\ell_2(\mathbb N_0,\mathbb C)$. The Arveson space corresponds to the restriction of
  $z$ and $w$ to the unit ball of $\mathbb C^N$.\\
  
To verify that a given kernel is indeed a CNP-kernel we will use the following theorem by Th. Kaluza~\cite[Satz 3]{MR1544949};
see also \cite[Theorem 1.3 p. 2]{MR2869264} and \cite[Lemma 7.38 p. 90]{MR1882259}. We note that the motivation for the result in Kaluza's paper consists in the study of the coefficients of reciprocal power series.

\begin{theorem}\label{Kaluza}
Let $f(z)=\sum_{n=0}^\infty a_nz^n$ such that 
\begin{enumerate}[i)]
\item $a_0=1$ and $a_n>0$, $n=1,2,\ldots$,
\item $\frac{a_n}{a_{n-1}}\leq\frac{a_{n+1}}{a_n},$
\end{enumerate}
then the series $1/f(z)=1-\sum_{n=1}^\infty b_nz^n$ satisfies $b_n\geq 0$ for all $n=1,2,\dots$. 
\end{theorem}

Unfortunately, not all kernels satisfy this condition, e.g., we consider the case
\[
  a_n=\frac{1}{[n]_q!}
\]
Then,
\[
  \begin{split}
    \frac{a_n}{a_{n-1}}&=\frac{[(n-1)]_q! }{[n]_q!}\\
&=\frac{1}{1+q+\cdots +q^{n-1}}\\
&\not\le    \frac{a_{n+1}}{a_{n}},
    \end{split}
  \]
  since the latter equals $\frac{1}{1+q+\cdots +q^{n}}$. In fact the condition means that the sequence $(a_{n+1}/a_n)_n$ cannot decrease faster than a geometric sequence. 

A standard example of a CNP-kernel is the Dirichlet kernel
  \[
    K(z,w)=\sum_{n=0}^\infty\frac{z^n\overline{w}^n}{n+1}.
  \]
Here we have
  \[
a_n^2=\frac{1}{(n+1)^2}\le a_{n-1}a_{n+1}=\frac{1}{n(n+2)},
\]
since $n(n+2)\le (n+1)^2$. This example appears in \cite[Satz 2 p. 162]{MR1544949}. It is an important example since the
number of negative power coefficients in its inverse  (i.e. number of negative squares of the corresponding kernel) is infinite.

Another interesting example is the Hardy-Sobolev space of functions over the unit disk where one imposes that not only the function on the boundary is in $L^2$, but also its derivatives. These function spaces have been used in \cite{BZ} to study the recovery of analytic function from its boundary values.
Classically, they are reproducing kernel Hilbert spaces with the reproducing kernel
\[
K_n(z,w)=\sum_{k=0}^\infty \frac{z^k\overline{w^k}}{1+k^2+k^2(k-1)^2+\cdots+k^2\cdots (k-n+1)^2}.
\]
In particular, for the case $n=1$ we have that the Hardy-Sobolev space where $f$ and $f^\prime$ is in $H_2$ has the reproducing kernel
  \[
K_1(z,w)=\sum_{n=0}^\infty \frac{z^n\overline{w^n}}{1+n^2}
    \]
which is not an CNP kernel. Nevertheless, we can consider the modified kernel
    \[
K_{1,\epsilon}(z,w)=\sum_{n=0}^\infty \frac{z^n\overline{w^n}}{1+\frac{n^2}{\e}}
  \]
which is a CNP kernel for $\e\le 1/2$.

This can be easily seen from the following calculation:
  \[
    \frac{1}{(1+\frac{n^2}{\e})^2}\le \frac{1}{(1+\frac{(n-1)^2}{\e})}\frac{1}{(1+\frac{(n+1)^2}{\e})}
  \]

  \[
   (1+\frac{n^2}{\e})^2\ge (1+\frac{(n-1)^2}{\e})(1+\frac{(n+1)^2}{\e})
 \]
 \[
  1+\frac{2n^2}{\e}+\frac{n^4}{\e^2}\ge1+\frac{2n^2+2}{\e}+\frac{(n^2-1)^2}{\e^2}\\
\]
\[
0\ge \frac{2}{\e}+\frac{1-2n^2}{\e^2}
\]
which leads to
\[
2\e+1\le 2n^2.
\]
Therefore, using Theorem~\ref{Kaluza} we have that the kernel is CNP for $\e\le 1/2$.

  \begin{lemma}
Let $(a_n)$ and $(b_n)$ be two sequences satisfying the conditions of Theorem \ref{Kaluza}, and let $c>0$. The the sequences $(a_n^c) $ and $(a_nb_n)$ satisfy also the conditions of the theorem. 
\end{lemma}

The sequence $(a_nb_n)$ is called the Hadamard product of the sequences; see \cite{hadamard}. This example, and more involved ones, appear in
\cite[Theorem 2.11 p. 9]{MR2869264}.\\

Hence the functions
\[
K(z,w)=\sum_{n=0}^\infty\frac{z^n\overline{w}^n}{(n+1)^2} \quad{\rm and}\quad K(z,w)=\sum_{n=0}^\infty\frac{z^n\overline{w}^n}{\sqrt{n+1}} 
  \]
are CNP kernels.\\

We can even state the following general theorem.

\begin{theorem}
Let $f:\mathbb{R}_+\to\mathbb{R}$ a continuously twice differentiable positive function which satisfies
\begin{equation}
(f^\prime(x))^2\leq f^{\prime\prime}(x)f(x),\quad x\in\mathbb R_+.
  \end{equation}
then the kernel $K(z,\overline{w})=\sum_{n=0}^\infty a_n z^n$ with $a_n=f(n)$ is a CPN-kernel.
\end{theorem}

\begin{proof}
To use Theorem~\ref{Kaluza}  we need to verify the condition
$$
\frac{f(x)}{f(x-1)}\leq \frac{f(x+1)}{f(x)}
$$
or with other words, the function $\frac{f(x+1)}{f(x)}$ is monotone increasing. This means that we have for the derivative the condition
\[
f^\prime(x+1)f(x)-f(x+1)f^\prime(x)\ge 0.
\]
Thus, the function $f^\prime/f$ has to be monotone increasing. Now, 
$$
\left(\frac{f^\prime}{f}\right)^\prime=\frac{f^{\prime\prime}f-(f^\prime)^2}{f}
$$ 
has to be positive. This means that $(f^\prime)^2\leq f^{\prime\prime}f$.

\end{proof}

We can even give the following interesting example by considering the partition function from statistical physics
$$
f(x)=\sum_k e^{x E_k}
$$
with given constants $E_k\geq 0$. Then the condition $\frac{f(n)}{f(n-1)}\leq \frac{f(n+1)}{f(n)}$
takes the form
$$
\frac{\sum_k e^{n E_k}}{\sum_k e^{n E_k}e^{-E_k}}\leq \frac{\sum_k e^{n E_k}e^{E_k}}{\sum_k e^{n E_k}}
$$
which leads to 
$$
\left(\sum_k e^{n E_k}\right)^2\leq \left({\sum_k e^{n E_k}e^{-E_k}}\right)\left({\sum_k e^{n E_k}e^{E_k}}\right)
$$
wich can be rewritten as
$$
\left(\sum_k e^{n E_k}\right)^2=\left(\sum_k e^{n \frac{E_k}{2}}e^{\frac{-E_k}{2}} e^{\frac{E_k}{2}} e^{n \frac{E_k}{2}}\right)^2\leq \left({\sum_k e^{n E_k}e^{-E_k}}\right)\left({\sum_k e^{n E_k}e^{E_k}}\right),
$$
but this is just Cauchy-Schwarz inequality. 

\section{CNP kernels in $q$-analysis}
\setcounter{equation}{0}

Let us now study the $q$-counterparts of various CNP-kernels.
We start with the  $q$-Dirichlet kernel. By replacing $n$ in the denominator of the running term of the Dirichlet kernel  by $[n]_q=\frac{1-q^n}{1-q}$ we have (after dividing by $1-q$)
  \[
K_q(z,w)=\sum_{n=1}^\infty\frac{z^n\overline{w}^n}{1-q^n}
    \]
  which can be rewritten as
    \[
K_q(z,w)=\sum_{a=0}^\infty\frac{z\overline{w}q^a}{1-z\overline{w}q^a}.
\]
We check Theorem \ref{Kaluza} with $a_n=\frac{1}{1-q^n}$ and get:
\[
\begin{split}
  \dfrac{\dfrac{1}{1-q^n}}{\dfrac{1}{1-q^{n-1}}}&\le \dfrac{\dfrac{1}{1-q^{n+1}}}{\dfrac{1}{1-q^{n}}}\\
        &\iff\\
        \frac{1-q^{n-1}}{1-q^n}&\le         \frac{1-q^{n}}{1-q^{n+1}}\\
        &\iff\\
        (1-q^{n-1})(1=q^{n+1})&\le (1-q^n)^2\\
        &\iff\\
        2&\le \frac{1}{q}+q
\end{split}
\]
which always hold.\\

Another example is the following $q$-version of the kernel studied in~\cite{ak4}. This kernel is given by 
$$
K_{q,r}(z,w)=\sum_{n=0}^\infty \frac{\Gamma_q(n+r)}{\Gamma_q(r)[n]_q!}z^n\overline{w}^n
$$
where $\Gamma_q$ denotes the $q$-Gamma function
$$
\Gamma_q(x)=(1-q)^{1-x}\Pi_{n=0}^\infty\frac{1-q^{n+1}}{1-q^{n+x}}=(1-q)^{1-x}\frac{(q;q)_\infty}{(q^x,q)_\infty}, \qquad |q|<1
$$
and $0\leq r\leq 1$.

Using Theorem~\ref{Kaluza} we can check if the above kernels are CNP-kernels. With $a_n=\frac{\Gamma_q(n+r)}{\Gamma_q(r)[n]_q!}$ the condition
$$
\frac{a_n}{a_{n-1}}\leq\frac{a_{n+1}}{a_n}
$$
has to be satisfied. From
$$
\frac{a_n}{a_{n-1}}\leq\frac{a_{n+1}}{a_n}\Leftrightarrow \frac{\frac{\Gamma_q(n+r)}{\Gamma_q(r)[n]_q!}}{\frac{\Gamma_q(n+r-1)}{\Gamma_q(r)[n-1]_q!}}\leq\frac{\frac{\Gamma_q(n+r+1)}{\Gamma_q(r)[n+1]_q!}}{\frac{\Gamma_q(n+r)}{\Gamma_q(r)[n]_q!}}
$$
we have
$$
1-q^{n+r+1}-q^n-q^{2n+r+1}\geq 1-q^{n+1}-q^{n+r}+q^{2n+r+1}
$$
which is equivalent to 
$$
q^n(1+q^r-q-q^{r-1})=q^n(1-q)(1-q^{r})\leq 0.
$$
Since the last inequality is obviously wrong for $q\neq 0$ our condition is not satisfied for $q\neq 0$ which is the case studied in~\cite{ak4}.

Furthermore, we also want to consider a q-version of Hardy-Sobolev spaces. 
Here we consider spaces with the reproducing kernel
$$
K_{n, \epsilon}(z,w)=\sum_{k=0}^\infty \frac{z^k\overline{w}^k}{(1-q)+([k]_q)^2+([k]_q)^2([k-1]_q)^2+\ldots+([k]_q)^2\cdots([k-n+1]_q)^2}, \qquad n\geq 1.
$$
We remark that in difference to the classic Hardy-Sobolev spaces in the above definition we assume that $f$ and $R_qf$ belong to the $q$-Fock space $\mathcal{H}_{2,q}$.  Of course, in the case $q=0$ we do not get the usual Hardy-Sobolev spaces, but rather spaces where $f$ and the backward shift-operator $R_0f$ belong to the Hardy space.  

Now, for $n=1$ we have the modified kernel
$$
K_{1, \epsilon}(z,w)=\sum_{k=0}^\infty \frac{z^k\overline{w}^k}{(1-q)+\frac{([k]_q)^2}{\epsilon}}.
$$

Here we get that
  \[
    \frac{1}{(1+\frac{[n]_q^2}{\e})^2}\le \frac{1}{(1+\frac{[n-1]_q^2}{\e})}\frac{1}{(1+\frac{[n+1]^2}{\e})}
  \]
which means
  \[
   (\epsilon+[n]_q^2)^2\ge (\epsilon+[n-1]_q^2)(\epsilon+[n+1]_q^2)
 \]
or
  \[
   2\epsilon [n]_q^2+[n]_q^4\ge \epsilon ( [n-1]_q^2+[n+1]_q^2)+[n-1]_q^2 [n+1]_q^2
 \]
 
 $$\epsilon ( 2  [n]_q^2 - [n-1]_q^2 - [n+1]_q^2) \geq   [n-1]_q^2 [n+1]_q^2 - [n]_q^4
$$

~

$$\epsilon \left(  \frac{2(1-q^{n})^2 - (1-q^{n-1})^2 - (1-q^{n+1})^2}{(1-q)^2}\right) \geq  \frac{(1-q^{n-1})^2 (1-q^{n+1})^2 - (1-q^{n})^4}{(1-q)^4}$$

Note that
\begin{gather*}
(1-q^{n})^2 - (1-q^{n-1})^2 = (1-q^{n} -1+q^{n-1})(1-q^{n} +1-q^{n-1}) \\
= q^{n-1}(1-q) (2 - q^{n-1} (1+q))
\end{gather*}
so that the left-hand-side is 
\begin{gather*}
\epsilon ( 2  [n]_q^2 - [n-1]_q^2 - [n+1]_q^2) = \epsilon \frac{q^{n-1}(1-q) (2 - q^{n-1} (1+q)) - q^{n}(1-q) (2 - q^{n} (1+q))}{(1-q)^2} \\
= \epsilon \frac{q^{n-1} [(2 - q^{n-1} (1+q)) - q(2 - q^{n} (1+q))]}{(1-q)}\\
= \epsilon (q^{n-1} [(2 - q^{n-1} ) - q(2 - q^{n} )]\frac{(1+q)}{(1-q)}\\
= \epsilon q^{n-1} (2 -2q - q^{n-1}  - q^{n+1} )\frac{(1+q)}{(1-q)}
\end{gather*}
Same for 
$$(1-q^{n-1})^2 (1-q^{n+1})^2 = (1 -q^{n-1} -q^{n+1} +q^{2n})^2$$
and the right hand side is
\begin{gather*}
\frac{(1-q^{n-1})^2 (1-q^{n+1})^2 - (1-q^{n})^4}{(1-q)^4} = \frac{(1 -q^{n-1} -q^{n+1} +q^{2n})^2 - (1 -2q^{n} +q^{2n})^2}{(1-q)^4} \\
= \frac{(1 -q^{n-1} -q^{n+1} +q^{2n}  - 1 +2q^{n} -q^{2n})(1 -q^{n-1} -q^{n+1} +q^{2n} + 1 -2q^{n} +q^{2n})}{(1-q)^4} \\
= \frac{( 2q^n - q^{n-1}-q^{n+1} )(2+2q^{2n} -q^{n-1}  - 2q^{n} -q^{n+1}  )}{(1-q)^4} \\
= -q^{n-1} \frac{( 1-2q +q^{2} )(2+2q^{2n} -q^{n-1}  - 2q^{n} -q^{n+1}  )}{(1-q)^4} \\
= -q^{n-1} \frac{(2+2q^{2n} -q^{n-1}  - 2q^{n} -q^{n+1}  )}{(1-q)^2} 
\end{gather*}

Combining we have 
\begin{gather*}
\epsilon q^{n-1}  (2 -2q - q^{n-1}  - q^{n+1} )\frac{(1+q)}{(1-q)} \geq -  q^{n-1} \frac{(2+2q^{2n} -q^{n-1}  - 2q^{n} -q^{n+1}  )}{(1-q)^2} \\
\epsilon (2 -2q - q^{n-1}  - q^{n+1} )(1+q) \geq \frac{q^{n-1} (1+q)^2  -2(1
+q^{2n}) }{(1-q)} \geq q^{n-1} (1+q)^2  -2(1
+q^{2n})
\end{gather*} since $\frac{1}{1-q} > 1.$

Or use instead 
\begin{gather*}
\epsilon (2 -2q - q^{n-1}  - q^{n+1} )(1+q) \geq \frac{q^{n-1} (1+q)^2  -2(1
+q^{2n}) }{(1-q)} \\
\epsilon (2 -2q - q^{n-1}  - q^{n+1} )(1-q^2) \geq q^{n-1} (1+q)^2  -2(1
+q^{2n})  
\end{gather*}

This means that our condition for $\epsilon$ is 
$$
\epsilon (2 -2q - q^{n-1}  - q^{n+1} )(1-q^2)\geq \frac{q^{n-1} (1+q)^2  -2(1
+q^{2n})} {2 -2q - q^{n-1}  - q^{n+1} }.
$$

To see study the condition we first look at the sign of the nominator and denominator of the left-hand side. 
For the denominator we get that it is negative if 
$$
2 -2q - q^{n-1}  - q^{n+1}\leq 0 
$$
the
$$
2(1-q)\leq	q^n\left(\frac{(1-q)(1+q)}{q}\right) \Leftrightarrow 2\leq q^{n-1}(1+q)
$$
which leads to the condition
$$
n\leq 1-\mathrm{log}_{1/q}\left(\frac{2}{1+q}\right).
$$
For the nominator to be negative we have the condition
$$
q^{n-1} (1+q)^2  -2(1+q^{2n})\leq 0 \Leftrightarrow q^{n-1}(1+q)^2\leq 2(1+q^{2n})
$$
which can be rewritten as 
$$
\frac{q^{n+1}+q^{n-1}}{1+q^{2n}}\leq 2.
$$
With $\alpha=q^{n+1}$ and $\beta=q^{n-1}$ we have 
$$
\frac{\alpha+\beta}{1+\alpha\beta}\leq 1\leq 2,
$$
since $\alpha<1$ and $\beta<1$.



\section{Interpolation for Schur multipliers associated CNP kernels}
\setcounter{equation}{0}
This last section is essentially of a survey type. The first remark is that, given a CNP kernel of the form \eqref{cnp00},
\[
\frac{1}{1-\langle b(z),b(w)\rangle_{\ell_2(A,\mathbb C)}},
\]
a (say matrix-valued) function $S$ is a
Schur multiplier if and only if it is of the form $S(z)=s(p(z))$ where $s$ is a Schur multiplier of the associated Arveson space. The interpolation problem for a CNP kernel as stated below in Problem~\ref{111} reduces thus to an 
interpolation problem in the Arveson space. As in the case of a single complex variable, interpolation in the class of Schur multipliers of the Arveson
space can be performed using various methods. In this section we refer to the paper of Ball, Trent and Vinnikov \cite{btv}, who use the lurking
isometry method, and to the paper \cite{abk}, which uses the reproducing kernel method. The emphasis in this last paper was the Schur algorithm and
the interpolation problems considered there consists of the tangential Nevanlinna-Pick interpolation problem. We refer also to \cite{alpay2023complete} for
related discussions and references.\smallskip

Following \cite{abk} we discuss the solution of the following interpolation problem:

\begin{problem}
  \label{111}
Given a CNP kernel of the form $k(z,w)=\frac{1}{1-\langle p(z),p(w)\rangle}$, where $p$ is a map from some set $\Omega$ into the open unit ball of $\mathbb C^N$  (or of $\ell^2(\mathbb N,\mathbb C)$ if $N=\infty$) and given $N$ triples
$(w^{(1)},\xi_1,\eta_1),\ldots, (w^{(m)},\xi_m,\eta_m)$ in $\Omega\times\mathbb C^{p}\times \mathbb C^q$, describe the set of all (if any) Schur multipliers associated to $k(z,w)$ and such that
\begin{equation}
  S(w_j)^*\xi_j=\eta_j,\quad j=1,\ldots, m,
\end{equation}
hold.
\end{problem}

Assuming the matrix
\[
G=\left(\frac{I_p-S_jS_\ell^*}{1-\langle p(w_j),p(w_\ell)\rangle}\right)_{j,\ell=1}^m>0,
\]
set
\[
  C=\begin{pmatrix}\xi_1^*&\xi_2^*&\cdots&\xi_m^*\\
    \eta_1^*&\eta_2^*&\cdots&\eta_m^*\end{pmatrix},\qquad \mathbf J=\begin{pmatrix}I_{nN}&0\\0&J\end{pmatrix},
\]
and
\[
A_j={\rm \,diag}(\overline{w_j^{(1)}},\ldots, \overline{w_j^{(N)}}),\quad j=1,\ldots, m
  \]
\[
\Theta(\lambda)=\begin{pmatrix}0&I_{p+q}\end{pmatrix}+C\left(I_N-\sum_{j=1}^m\lambda_jA_j^*\right)G^{-1}\begin{pmatrix}(\lambda_1I_n-A_1^*)G^{1/2}&\cdots&\lambda_1I_n-A_m^*)G^{1/2}&-C^*J\end{pmatrix}
\]
Then,
\begin{equation}
  \frac{J-\Theta(\lambda){\mathbf J}\Theta(\mu)^*}{1-\langle \lambda, \mu\rangle}=
  C\left(I_N-\sum_{j=1}^m\lambda_jA_j\right)^{-1}G^{-1}\left(I_N-\sum_{j=1}^m\overline{\mu_j}A_j^*\right)^{-1}.
  \end{equation}
Furthermore
with
\[
\Theta=(\theta_{ij})_{i,j=1}^2
\]
where $\theta_{11}$ is $\mathbb C^{p\times p}$-valued, the linear fractional transformation
\begin{equation}
S(\lambda)=\left(\theta_{11}(p(\lambda))\sigma(\lambda)+\sigma_{12}(\lambda)\right)\left(\theta_{21}(p(\lambda))\sigma(\lambda)+\sigma_{22}(\lambda)\right)^{-1}
\end{equation}
defines the set of all solutions of Problem \ref{111} when $\sigma$ varies along the $\mathbb C^{p\times (q+N-1)}$-valued Schur multipliers of
$(\mathcal H(k))^{p\times p}$. See \cite[\S8]{abk} for these. Note that the formulas are not optimal in the sense that for $N=1$, the matrix $\mathbf J$
does not reduce to $J$.

\section*{Acknowledgments} It is a pleasure to thank Prof. Vladimir Bolotnikov for mentioning to us Lemma \ref{vb}.\smallskip

{\bf Funding:} D. Alpay thanks the Foster G. and Mary McGaw Professorship in
  Mathematical Sciences, which supported his research. The second and third author were supported by Portuguese funds through the CIDMA - Center for Research and Development in Mathematics and Applications, and the Portuguese Foundation for Science and Technology (``FCT--Funda\c{c}\~ao para a Ci\^encia e a Tecnologia''), within project UIDB/04106/2020 and UIDP/04106/2020.\smallskip

{\bf Data Availability} Not applicable.\smallskip

{\bf Declarations Conflict of interest} The authors have no other relevant financial or non-financial interests to disclose.
\bibliographystyle{plain}
  \def\lfhook#1{\setbox0=\hbox{#1}{\ooalign{\hidewidth
  \lower1.5ex\hbox{'}\hidewidth\crcr\unhbox0}}} \def\cprime{$'$}
  \def\cprime{$'$} \def\cprime{$'$} \def\cprime{$'$} \def\cprime{$'$}
  \def\cprime{$'$}

\end{document}